\documentclass{amsart}

\usepackage{amssymb}
\usepackage{amsfonts}
\usepackage{amsmath}
\usepackage{graphicx}
\usepackage{mathrsfs}
\usepackage{dsfont}
\usepackage{amscd}
\usepackage{multirow}
\usepackage[all]{xy}
\normalfont
\usepackage[T1]{fontenc}
\usepackage{verbatim}
\usepackage{hyperref}

\newcommand{\N}{\mathds{N}}

\newcommand{\R}{\mathds{R}}
\newcommand{\C}{\mathds{C}}
\newcommand{\T}{\mathds{T}}

\renewcommand{\H}{\mathscr{H}}
\newcommand{\F}{\mathds{F}}
\newcommand{\id}{\mathds{1}}
\newcommand{\A}{\mathscr{A}}
\newcommand{\Iso}{\mathscr{I}}

\theoremstyle{plain}
\newtheorem{theorem}{Theorem}[section]

\theoremstyle{definition}
\newtheorem{definition}[theorem]{Definition}

\theoremstyle{remark}

\numberwithin{equation}{section}

\begin{document}
\title[Classification of Noncommutative domains]{Classification of Noncommutative Domain Algebras}
\author{Alvaro Arias}
\address{Department of Mathematics \\ University of Denver \\ Denver CO 80208}
\email{aarias@math.du.edu}
\urladdr{http://www.math.du.edu/\symbol{126}aarias}

\author{Fr\'{e}d\'{e}ric Latr\'{e}moli\`{e}re}
\address{Department of Mathematics \\ University of Denver \\ Denver CO 80208}
\email{frederic@math.du.edu}
\urladdr{http://www.math.du.edu/\symbol{126}frederic}

\date{\today}

\begin{abstract}
Noncommutative domain algebras are noncommutative analogues of the algebras of holomorphic functions on domains of $\C^n$ defined by holomorphic polynomials, and they generalize the noncommutative Hardy algebras. We present here a complete classification of these algebras based upon techniques inspired by multivariate complex analysis, and more specifically the classification of domains in hermitian spaces up to biholomorphic equivalence.
\end{abstract}
\maketitle

\section{Introduction}

Noncommutative domain algebras, introduced by Popescu in \cite{Popescu08} are universal objects in the category of operator algebras for certain polynomial relations and are noncommutative analogues of algebras of holomorphic functions on domains in hermitian spaces, as well as generalizations of Hardy Algebras. This note presents a classification of these algebras up to completely isometric isomorphism, based upon their defining symbols, and thus complete the work initiated by the authors in \cite{Latremoliere09b} and \cite{Latremoliere11b}. Our methods are based on the deep interplay between analysis of this type of operator algebras and analysis of multivariate holomorphic functions. This note also exploits an observation made in \cite{Davidson11} which helps us conclude Theorem \ref{main}.

Algebras of weighted shifts on Fock spaces \cite{Popescu91,Arveson98,Pitts98,Davidson01} find their origin in the study of row contractions, dilations, and commutant lifting theorems. Their applications became important in interpolation problems \cite{Davidson98, Arias00} and the study of $E_0$-semigroups \cite{Davidson11} or noncommutative complex analysis \cite{Latremoliere09b,Latremoliere10,Latremoliere11b}. 

In \cite{Latremoliere09b} and \cite{Latremoliere11b}, the authors of this note studied the geometry of the spectra of noncommutative domains. They used Thullen's \cite{Thullen31} and Sunada's \cite{Sunada78} classification of Reinhardt domains and combinatorial arguments to fully classify a large class of domain algebras. However, one important case remained unsolved, and this note present the solution for this last unresolved class, thus completing the classification of all noncommutative domain algebras with polynomial symbols.

Noncommutative domain algebras are naturally --- and best --- represented as norm closed algebras of operators on full Fock spaces. For any (complex) Hilbert space $\H$, the associated full Fock space, denoted by $\mathscr{F}(\H)$, is the Hilbert sum $\oplus_{n\in \N} \H^{\otimes n}$ with the convention that $\H^0 = \C$. It is useful for the definition of noncommutative domain algebras to introduce the following alternative description of $\mathscr{F}(\H)$. Let $\F_+^n$ be the free semigroup over $\{1,\ldots,n\}$ whose neutral element will be denoted by $\emptyset$ by abuse of notation. For any elements $t_1,\ldots,t_n$ of any unital associative algebra with unit $\id$, and for any nonempty word $\alpha = i_1\cdots i_k$ in $\F_+^n$, with $i_1,\ldots,i_k \in \{1,\ldots,k\}$, we denote by $t_\alpha$ the product $\prod_{j=1}^k t_{i_j}$, while $t_\emptyset$ will be set to $\id$ by convention. In particular, $\mathscr{F}(\C^n)$ is a unital associative algebra for the tensor product. We fix henceforth the canonical basis $\{e_1,\ldots,e_n\}$ of $\C^n$ and observe that $(e_\alpha)_{\alpha\in\F_+^n}$ is a Hilbert basis of $\mathscr{F}(\C^n)$, which allows us to identify $\mathscr{F}(\C^n)$ with $\ell^2(\F_+^n)$.

Fix $n\in\N,n>0$. Let $\mathscr{P}_n = \C[X_1,\ldots,X_n]$ be the algebra of polynomials in $n$ noncommuting indeterminates $X_1,\ldots,X_n$ --- which is the free associative algebra generated by $\F_+^n$, where each generator $i\in \{1,\ldots,n\}$ of $\F_+^n$ is identified with $X_i$. Thus any element of $\mathscr{P}_n$ is of the form $\sum_{\alpha\in\F_+^n}a_\alpha X_\alpha$ for an almost zero family $(a_\alpha)_{\alpha\in\F_+^n}$ of complex numbers. For any $f=\sum_{\alpha\in\F_+^n}a_\alpha X_\alpha \in \mathscr{P}_n$ and any Hilbert space $\H$ whose C*-algebra of bounded linear operators will be denoted by $\mathscr{B}(\H)$, the noncommutative domain $\mathscr{D}_f(\H)$ is defined by:
\[
\mathscr{D}_f(\H) = \left\{ (T_1,\ldots,T_n) \in \mathscr{B}(\H)^n : \sum_{\alpha\in\F_+^n} a_\alpha T_\alpha T_\alpha^\ast \leq \id_\H \right\}\text{,}
\]
where we used the notation $T_\alpha^\ast$ for $(T_\alpha)^\ast$. 

If $f$ is an element of the subset $\mathscr{S}_n\subseteq \mathscr{P}_n$ defined by:
\[
\mathscr{S}_n = \left\{ f=\sum_{\alpha\in\F_+^n}a_\alpha X_\alpha \left|
a_\emptyset = 0, \forall \alpha \in \F_+^n \;\; a_\alpha \geq 0\wedge \left(|\alpha|=1\implies a_\alpha>0\right)\right.
  \right\} \text{,}
\]
where $|\alpha|$ is the word-length of $\alpha\in\F_+^n$, one can construct an explicit, universal $n$-tuples of operators in $\mathscr{D}_f(\mathscr{F}(\C^n))$. We shall refer to elements of $\mathscr{S}_n$ as \emph{$n$-symbols}, while Popescu refers to them as \emph{positive regular $n$-free formal polynomials}. Let us fix an $n$-symbol $f=\sum_{\alpha\in\F_+^n}a_\alpha X_\alpha$. We define the weighted shifts $W_1^f,\ldots,W_n^f$ on $\mathscr{F}(\C^n)$  by extending:
\begin{equation}\label{W}
W_j^f (e_\alpha) = \sqrt{\frac{b_\alpha}{b_{j\cdot \alpha}}} (e_j\otimes e_{\alpha}) \text{\,where\;} b_\alpha= \sum_{k=1}^{|\alpha|}\sum_{\substack{\gamma_1\cdots\gamma_k = \alpha \\ |\gamma_1|\geq 1,\cdots,|\gamma_k|\geq 1}} a_{\gamma_1}\cdots a_{\gamma_k} \text{for all \;}\alpha \in \F_+^n\text{.}
\end{equation}
We then have the following fundamental universal property:
\begin{theorem}(Popescu, \cite{Popescu08})\label{U}
Let $n\in\N,n>0$ and let $f\in\mathscr{S}_n$. Let $\A_f$ be the norm closure of the associative algebra generated by the set of operators $\{W_1^f,\ldots,W_n^f\}$ defined by (\ref{W}) on $\mathscr{F}(\C^n)$. We have $(W_1^f,\ldots,W_n^f)\in\mathscr{D}_f(\mathscr{F}(\C^n))$. Moreover, for any Hilbert space $\H$ and any $(T_1,\ldots,T_n)\in\mathscr{D}_f(\H)$, there exists a unique completely contractive algebra morphism $\varphi$ from $\A_f$ onto the norm closure of the algebra generated by $T_1,\ldots,T_n$ such that $\varphi(W_j^f)=T_j$ for $j=1,\ldots,n$.
\end{theorem}

The purpose of this note is to completely classify the noncommutative domain algebras $\A_f$ defined for $f \in \mathscr{S}_n$ in Theorem (\ref{U}) in term of their symbol $f$.

\section{Main result}

We define the following equivalence relation on the set $\mathscr{S}_n$ of symbols \cite[Definition 2.1]{Latremoliere11b}:
\begin{definition}
Let $n\in\N,n>0$. Two elements $f,g \in \mathscr{S}_n$ are scale-permutation equivalent when there exists a permutation $\sigma$ of $\{1,\ldots,n\}$ and scalars $\lambda_1,\ldots,\lambda_n \in \R$ such that $f(X_1,\ldots,X_n) = g(\lambda_1 X_{\sigma(1)}, \ldots, \lambda_n X_{\sigma(n)} )$.
\end{definition}

Our main theorem is the following complete classification result:

\begin{theorem}\label{main}
Let $f\in \mathscr{S}_n$ and $g\in\mathscr{S}_m$ be two symbols. The noncommutative domain algebras $\A_f$ and $\A_g$ are completely isometrically isomorphic if and only if $n=m$ and $f$ and $g$ are scale-permutation equivalent.
\end{theorem}

Before we provide the proof of this result, let us recall from \cite{Latremoliere09b} the following fundamental duality construction upon which our work relies. Let $f\in\mathscr{S}_n$ for some $n\in\N,n>0$. Let $k\in\N,k>0$. Let $T=(T_1,\ldots,T_n)\in\mathscr{D}_f(\C^k)$. By universality of $\A_f$, there exists a completely contractive morphism $\left<T,\cdot\right>_k : \A_f \rightarrow \mathscr{B}(\C^k)$ such that $\left<T,W_j^f\right>_k = T_j$ for $j=1,\ldots,n$. Moreover, the function $\left<\cdot,a\right>_k$ is holomorphic on the interior of $\mathscr{D}_f(\C^k)$ and extends by continuity to $\mathscr{D}_f(\C^k)$. We shall abuse the terminology and call such a function holomorphic on the compact $\mathscr{D}_f(\C^k)$.

\begin{proof}
In this proof, isomorphisms are always meant for completely isometric isomorphisms of operator algebras.

By \cite[Lemma 4.4]{Latremoliere09b}, whenever two symbols are scale-permutation equivalent, their associated noncommutative domains are isomorphic. It is thus sufficient to prove the converse here.

Let us denote the set of all isomorphisms from $\A_f$ to $\A_g$ by $\Iso(\A_f,\A_g)$. By \cite[Theorem 3.7]{Latremoliere09b}, for any $\Psi\in\Iso(\A_f,\A_g)$, there exists a necessarily unique biholomorphic map $\widehat{\Psi}:\mathscr{D}_g(\C)\rightarrow \mathscr{D}_f(\C)$ such that, for all $\lambda=(\lambda_1,\cdots,\lambda_n)\in\mathscr{D}_g(\C)$ and $a\in\A_f$, we have $\left<\lambda,\Psi(a)\right>_1 = \left<\widehat{\Psi}(\lambda),a\right>_1 \text{.}$ Moreover, by \cite[Theorem 3.18]{Latremoliere09b} and \cite[Theorem 3.2]{Latremoliere11b}, if $\widehat{\Psi}(0)=0$ then $f$ and $g$ are scale-permutation equivalent.

Thus, let us assume that there exists an isomorphism $\Phi : \A_f \rightarrow \A_g$. In \cite[Theorem 3.4]{Latremoliere11b}, we were able to show that if either $\mathscr{D}_f(\C)$ or $\mathscr{D}_g(\C)$ is not biholomorphic to the unit ball $\mathds{B}_n = \{ (z_1,\ldots,z_n) \in \C^n : \sum_{j=1}^n z_j\overline{z_j} \leq 1\}$ of $\C^n$, then $\widehat{\Psi}(0)=0$, which implies that $f$ and $g$ are scale-permutation equivalent. Thus, it remains to consider the case where $\mathscr{D}_f(\C)$ and $\mathscr{D}_g(\C)$ are $\mathds{B}_n$, up to replacing $f$ and $g$ by scale-permutation equivalent symbols. 

Let $\omega = \widehat{\Phi}^{-1}(0)$. If $\omega=0$ then by \cite[Theorem 3.2]{Latremoliere11b}, we can already conclude that $f$ and $g$ are scale-permutation equivalent. Henceforth, we shall assume $\omega \not=0$. We adapt an argument in \cite{Davidson11}. Let $D_g = \C \omega \cap \mathds{B}_n$ and $D_f = \widehat{\Phi}(D_g)$. By construction, $D_g$ is a disk in the plane $\C \omega$. Now, as a conformal self-map of the unit ball, $\widehat{\Phi} = \Upsilon \circ \varphi_{\omega}$ where $\varphi_{\omega}$ is the conformal map $\varphi_{\omega} = \frac{\psi_{\omega}}{|\psi_{\omega}|^2}$ where $\psi_{\omega} : z \in \mathds{B}_n \mapsto \omega + (1-|\omega|^2)\frac{\omega-x}{|\omega-x|^2}$, and $\Upsilon$ is unitary. Indeed, $\varphi_{\omega}$ is a easily seen to be a conformal map such that $\varphi_{\omega}\circ\varphi_{\omega}$ is the identity and $\varphi_{\omega}(0)=\omega$. Thus $\widehat{\Phi}\circ\varphi_{\omega}(0) = 0$ and hence, by Cartan's Lemma, $\widehat{\Phi}\circ\varphi_{\omega}$ is a unitary $\Upsilon$. Hence $\widehat{\Phi} = \widehat{\Phi}\circ \varphi_{\omega} \circ \varphi_{\omega} = \Upsilon\circ \varphi_{\omega}$. Now, $\varphi_{\omega}$ maps the plane $\C\omega$ to itself by construction, so $D_f = \widehat{\Phi}(D_g) = \Upsilon(D_g)$, and since $\Upsilon$ is unitary, $D_f$ is also a disk.

We now set $G = \{ z \in D_g : \exists \Psi \in \Iso(\A_f,\A_g) \;\; \widehat{\Psi}(z)=0  \}$ and $F = \{ z \in D_f : \exists \Psi \in \Iso(\A_f,\A_f) \;\; \widehat{\Psi}(0)=z \}$. The first observation is that $G$ and $F$ are circular domains. Indeed, note first that for any $\lambda\in\C,|\lambda|=1$, we can extend the function which, to each $j\in \{1,\ldots,n\}$, maps $W_j^g$ to $\overline{\lambda} W_j^g$, to an automorphism $\Lambda$ of $\A_g$ by \cite[Lemma 4.4]{Latremoliere09b}. Now, for any $b\in G$ there exists $\Psi\in\Iso(\A_f,\A_g)$ such that $\widehat{\Psi}(b)=0$. We have $\widehat{\Lambda\circ \Psi} = \widehat{\Psi}\circ \widehat{\Lambda}$ so $\widehat{\Lambda\circ\Psi}(\lambda b)= 0$ and $\Psi\circ\Lambda \in \Iso(\A_f,\A_g)$, so $\lambda b \in G$ since $\lambda b \in D_g$ as $D_g$ is a disk, hence circular. The same argument of course applies to $F$.

Now, by construction, $\omega \in G$. Thus $G$ contains the circle $\T_g = \{ \lambda \omega : |\lambda| = 1\}$. Since $\widehat{\Phi}$ restricted to the disk $D_g$ is a M{\"o}bius map whose poles lie outside of $D_g$, it preserves circles. Let $\T_f = \widehat{\Phi}(\T_g)$. On the other hand, if $b \in G$ then there exists $\Psi\in\Iso(\A_f,\A_g)$ such that $\widehat{\Psi}(b)=0$, so $\widehat{\Phi}(b)=\widehat{\Phi}\circ\widehat{\Psi}^{-1}(0) = \widehat{\Psi^{-1}\circ\Phi}(0)$ and $\Psi\circ\Phi \in \Iso(\A_f,\A_f)$. Hence $\widehat{\Phi}(b) \in F$. Since $\omega\in \T_g$ we conclude that $0 \in \T_f$.

Hence, $F$ contains the circle $\T_f$ containing the origin. Since $F$ is circular, one checks easily that $F$ contains the interior of $\T_f$. Since $G=\widehat{\Phi}^{-1}(F)$ by as similar argument as above, we conclude that $G$ contains the interior of $\T_g = \widehat{\Phi}^{-1}(\T_f)$, which in turn contains $0$. Hence, there exists $\Psi \in \Iso(\A_f,\A_g)$ such that $\widehat{\Psi}(0)=0$. By \cite[Theorem 3.2]{Latremoliere11b}, we conclude that $f$ and $g$ are scale-permutation equivalent.
\end{proof}


\bibliographystyle{amsplain}
\bibliography{../thesis}

\end{document}